\documentclass[twoside]{amsart}
\usepackage{amsfonts,amsthm}
\usepackage{enumerate}
\theoremstyle{plain}
\newtheorem*{thm A}{Theorem~A}
\newtheorem*{thm B}{Theorem~B}
\newtheorem*{thm C}{Theorem~C}
\newtheorem*{thm D}{Theorem~D}
\newtheorem*{thm E}{Theorem~E}
\newtheorem*{Main Theorem}{Main Theorem}

\newtheorem*{thm 1}{Theorem~1}
\newtheorem*{thm 2}{Theorem~2}
\newtheorem*{pro A}{Proposition~A}
\newtheorem*{pro B}{Proposition~B}
\newtheorem*{lem A}{Lemma~A}
\newtheorem*{lem B}{Lemma~B}
\newtheorem*{lem C}{Lemma~C}
\newtheorem*{lem D}{Lemma~D}

\newtheorem*{proof*}{\it The proof of Corollary}
\newtheorem*{rem}{Remark}

\newtheorem{theorem}{Theorem}[section]
\newtheorem{corollary}[theorem]{Corollary}
\newtheorem{lemma}[theorem]{Lemma}
\newtheorem{proposition}[theorem]{Proposition}

\theoremstyle{plain}
\newcommand{\be}{\begin{equation}}
\newcommand{\ee}{\end{equation}}
\newcommand{\bea}{\begin{eqnarray}}
\newcommand{\eea}{\end{eqnarray}}
\newcommand{\ba}{\begin{array}}

\newcommand{\ea}{\end{array}}

\newcommand{\bc}{\begin{center}}
\newcommand{\ec}{\end{center}}

\newcommand{\benu}{\begin{enumerate}}
\newcommand{\eenu}{\end{enumerate}}
\newcommand{\bpr}{\begin{proposition}}
\newcommand{\epr}{\end{proposition}}
\newcommand{\ble}{\begin{lemma}}
\newcommand{\ele}{\end{lemma}}
\newcommand{\bco}{\begin{corollary}}
\newcommand{\eco}{\end{corollary}}

\def \Cal{\mathcal}

\def \ENK{{\eta}_{\nu}({\xi})}

\def \PKN{{\phi}{\xi}_{\nu}}

\def \NBo{SU_{2,m-1}/S(U_2{\cdot}U_{m-1})}
\def \NBt{SU_{2,m}/S(U_{2}{\cdot}U_{m})}
\def \CC{{\Bbb C}}
\def \RR{{\Bbb R}}
\def \HH{{\Bbb H}}

\def \ENK{{\eta}_{\nu}({\xi})}

\def \PKN{{\phi}{\xi}_{\nu}}

\def \NBo{SU_{2,m-1}/S(U_2U_{m-1})}
\def \NBt{SU_{2,m}/S(U_2U_m)}
\def \CC{{\Bbb C}}
\def \RR{{\Bbb R}}
\def \HH{{\Bbb H}}

\def \ENK{{\eta}_{\nu}({\xi})}

\def \PKN{{\phi}{\xi}_{\nu}}

\def \NBo{SU_{2,m-1}/S(U_2U_{m-1})}
\def \NBt{SU_{2,m}/S(U_2U_m)}
\def \CC{{\Bbb C}}
\def \RR{{\Bbb R}}
\def \HH{{\Bbb H}}

\begin{document}

\title[hypersurfaces in complex hyperbolic Grassmannians]
{Real HYPERSURFACES IN COMPLEX hyperbolic TWO-PLANE\\ GRASSMANNIANS
Related to the Reeb vector Field}
\vspace{0.2in}
\author[Young Jin Suh]{Young Jin Suh}
\address{\newline
Young Jin Suh
\newline Kyungpook National University,
\newline Department of Mathematics,
\newline Taegu 702-701, Korea}
\email{yjsuh@knu.ac.kr}


\thanks{* This work was supported by grants Proj. No.
NRF-2011-220-C00002 and Proj. No. NRF-2012-R1A2A2A-01043023 from National Research
Foundation of Korea.}

\begin{abstract}
 In this paper we give a characterization of real
hypersurfaces in noncompact complex two-plane Grassmannian
$SU_{2,m}/S(U_2U_m)$, $m \geq 2$ with Reeb vector field $\xi$
belonging to the maximal quaternionic subbundle $\Cal Q$. Then it becomes a tube over a
totally real totally geodesic ${\Bbb H}H^n$, $m=2n$, in noncompact
complex two-plane Grassmannian $SU_{2,m}/S(U_2U_m)$, a horosphere
whose center at the infinity is singular or another exceptional
case.
\end{abstract}

\maketitle

\section*{Introduction}
\setcounter{equation}{0}
\renewcommand{\theequation}{0.\arabic{equation}}
\vspace{0.13in}

Let us denote by $SU_{2,m}$ the set of
$(m+2){\times}(m+2)$-indefinite special unitary matrices
 and $U_m$ the set of $m{\times}m$-unitary matrices. Then
the Riemannian symmetric space $SU_{2,m}/S(U_2U_m)$, $m \geq 2$,
which consists of positive definite complex two-planes in indefinite
complex Euclidean space ${\Bbb C}_2^{m+2}$ (See page 315, Besse \cite{Be}), has a remarkable
feature that it is a Hermitian symmetric space as well as a
quaternionic K\"{a}hler symmetric space. In fact, among all
Riemannian symmetric spaces of noncompact type the symmetric spaces
$SU_{2,m}/S(U_2U_m)$, $m \geq 2$, are the only ones which are
Hermitian symmetric and quaternionic K\"{a}hler symmetric. So we will say such a Hermitian symmetric space of noncompact type
$SU_{2,m}/S(U_2U_m)$ a complex hyperbolic two-plane Grassmannian.
\par
\vskip 6pt
The existence of these two structures leads to a number of interesting
geometric problems on $SU_{2,m}/S(U_2U_m)$, one of which we are
going to study in this article. To describe this problem, we denote
by $J$ the K\"{a}hler structure and by ${\frak J}$ the quaternionic
K\"{a}hler structure on $SU_{2,m}/S(U_2U_m)$. Let $M$ be a connected
hypersurface in $SU_{2,m}/S(U_2U_m)$ and denote by $N$ a unit normal to $M$.
Then a structure vector field $\xi$ defined by ${\xi}=-JN$ is said to be a Reeb vector field.
\par
\vskip 6pt
Now let us denote by $TM$ the tangent
bundle of $M$. Then the maximal complex subbundle of $TM$ is defined by
${\Cal C} = \{ X \in TM \mid JX \in TM\}$, and the maximal
quaternionic subbundle ${\Cal Q}$ of $TM$ is defined by ${\Cal Q} =
\{ X \in TM \mid {\frak J}X \in TM\}$, where ${\frak J}={\text Span}\{J_1,J_2,J_3\}$, and $\{J_1,J_2,J_3\}$ denotes the quaternionic K\"ahler structure. The main subject we want to discuss in this paper is: {\it What can we say about real hypersurfaces in $SU_{2,m}/S(U_2U_m)$ with the Reeb vector field belonging to the maximal quaternionic subbundle ${\Cal Q}$}?
\par
\vskip 6pt
Before going to do this, we introduce some hypersurfaces in $SU_{2,m}/S(U_2U_m)$ with
invariant maximal complex subbundle $\Cal C$ and quaternionic subbundle $\Cal Q$ of $M$. Let us denote by $o \in SU_{2,m}/S(U_2U_m)$ the
unique fixed point of the action of the isotropy group $S(U_2U_m)$
on $SU_{2,m}/S(U_2U_m)$.
\par
\vskip 6pt
First we consider the conic (or geodesic) compactification of
$SU_{2,m}/S(U_2U_m)$. The points in the boundary of this
compactification correspond to equivalence classes of asymptotic
geodesics in $SU_{2,m}/S(U_2U_m)$. Every geodesic in this noncompact Grassmannian
lies in a maximal flat, that is, a two-dimensional Euclidean space embedded in $SU_{2,m}/S(U_2U_m)$ as
a totally geodesic submanifold. A geodesic in $SU_{2,m}/S(U_2U_m)$
is called singular if it lies in more than one maximal flat in
$SU_{2,m}/S(U_2U_m)$. A singular point at infinity is the
equivalence class of a singular geodesic in $SU_{2,m}/S(U_2U_m)$. Up
to isometry, there are exactly two singular points at infinity for
$SU_{2,m}/S(U_2U_m)$. The singular points at infinity correspond to
the geodesics in $SU_{2,m}/S(U_2U_m)$ which are determined by
nonzero tangent vectors $X$ with $JX \in {\frak J}X$
or $JX \perp {\frak J}X$ respectively.
\par
\vskip 6pt
Motivated by the results mentioned above, recently Berndt and the author \cite{BS3} have given a
complete characterization of horospheres in $SU_{2,m}/S(U_2U_m)$
whose center at infinity is singular as follows:
\par
\vskip 6pt

\begin{thm A}\label{Theorem A}
Let $M$ be a horosphere in
$SU_{2,m}/S(U_2U_m)$, $m \geq 2$. The following statements are
equivalent:
\par
(i) the center of $M$ is a singular point at infinity,
\par
(ii) the maximal complex subbundle ${\Cal C}$ of $TM$ is invariant
under the shape operator of $M$,
\par
(iii) the maximal quaternionic subbundle ${\Cal Q}$ of $TM$ is
invariant under the shape operator of $M$.
\end{thm A}
\vskip 6pt
\par
Next, we consider the standard embedding of $SU_{2,m-1}$ in
$SU_{2,m}$. Then the orbit $SU_{2,m-1} \cdot o$ of $SU_{2,m-1}$
through $o$ is the Riemannian symmetric space \hfill
\newline
$SU_{2,{m-1}}/S(U_2U_{m-1})$ embedded in $SU_{2,m}/S(U_2U_m)$ as a
totally geodesic submanifold. Every tube around
$SU_{2,m-1}/S(U_2U_{m-1})$ in $SU_{2,m}/S(U_2U_m)$ has the property
that both maximal complex subbundle ${\Cal C}$ and quaternionic subbundle ${\Cal Q}$ are invariant under the shape
operator.
\par
\vskip 6pt
Finally, let $m$ be even, say $m = 2n$, and consider the standard
embedding of $Sp_{1,n}$ in $SU_{2,2n}$. Then the orbit $Sp_{1,n}
\cdot o$ of $Sp_{1,n}$ through $o$ is the quaternionic hyperbolic
space ${\Bbb H}H^n$ embedded in $SU_{2,2n}/S(U_2U_{2n})$ as a
totally geodesic submanifold. Any tube around ${\Bbb H}H^n$ in
$SU_{2,2n}/S(U_2U_{2n})$ has the property that both ${\Cal C}$ and
${\Cal Q}$ are invariant under the shape operator.
\par
\vskip 6pt
As a converse of the statements mentioned above, we assert that
with one possible exceptional case there are no other such real
hypersurfaces. Related to such a result, we introduce another theorem due to Berndt and Suh \cite{BS3} as follows:
\begin{thm B}\label{Theorem B}
Let $M$ be a connected hypersurface in
$SU_{2,m}/S(U_2U_m)$, $m \geq 2$. Then the maximal complex subbundle
${\Cal C}$ of $TM$ and the maximal quaternionic subbundle ${\Cal Q}$
of $TM$ are both invariant under the shape operator of $M$ if and
only if $M$ is congruent to an open part of one of the following
hypersurfaces:
\par
$(A)$\quad a tube around a totally geodesic
$SU_{2,m-1}/S(U_2U_{m-1})$ in $SU_{2,m}/S(U_2U_m)$;
\par
$(B)$\quad a tube around a totally geodesic ${\Bbb H}H^n$ in
$SU_{2,2n}/S(U_2U_{2n})$, $m = 2n$;
\par
$(C)$\quad a horosphere in $SU_{2,m}/S(U_2U_m)$
 whose center at infinity is singular;
\par
or the following exceptional case holds:
\par
$(D)$\quad The normal bundle $\nu M$ of $M$ consists of singular
tangent vectors of type $JX \perp {\frak J}X$. Moreover, $M$ has at
least four distinct principal curvatures, three of which are given
by
\begin{equation*}
\alpha = \sqrt{2}\ ,\ \gamma = 0\ ,\ \lambda = \frac{1}{\sqrt{2}}
\end{equation*}
with corresponding principal curvature spaces
\begin{equation*}
T_\alpha = TM \ominus ({\Cal C} \cap {\Cal Q})\ ,\ T_\gamma = J(TM
\ominus {\Cal Q})\ ,\ T_\lambda \subset  {\Cal C} \cap {\Cal Q} \cap
J{\Cal Q}.
\end{equation*}
If $\mu$ is another (possibly nonconstant) principal curvature
function, then we have $T_\mu \subset {\Cal C} \cap {\Cal Q} \cap
J{\Cal Q}$, $JT_\mu \subset T_\lambda$ and ${\frak J}T_\mu \subset
T_\lambda $.
\end{thm B}
\par
\vskip 6pt
Usually, maximal complex subbundle of real hypersurfaces in a K\"ahler manifold is invariant under the shape operator
when the Reeb vector field $\xi =-JN$ is principal.
Accordingly, the maximal complex subbundle $\Cal C$ of $TM$ in Theorem B is
invariant under the shape operator if and only if the Reeb vector field
$\xi$ becomes a principal vector field for the shape operator $A$ of
$M$ in $SU_{2,m}/S(U_2U_m)$. In this case we call $M$ a Hopf hypersurface in $SU_{2,m}/S(U_2U_m)$.
\par
\vskip 6pt
Besides of this, a real hypersurface
$M$ in $\NBt$ also admits the maximal quatrnionic subbundle $\Cal Q$ and the orthogonal complement
${\Cal Q}^{\bot}$, which is spanned by {\it almost
contact 3-structure} vector fields $\{{\xi}_1,{\xi}_2,{\xi}_3\}$,
such that $T_xM={\Cal Q}{\oplus}{\Cal Q}^{\bot}$, $x{\in}M$.
\par
\vskip 6pt
In order to give some characterizations for hypersurfaces given in Theorem B, we \cite{S5} have considered a geometric condition that
the {\it Reeb flow on $M$ is isometric}, that is, the shape opeartor of $M$ in
$\NBt$ commutes with the structure tensor $\phi$. By virtue of this condition, we gave a characterization of real hypersurfaces of type $(A)$ or one of type
$(C)$ in $\NBt$. Historically, many geometers considered such a notion on several kinds of manifolds.  As a first the {\it isometric Reeb flow} on real hypersurfaces in complex projective space ${\Bbb C}P^n$ was investigated by Okumura \cite{O} , and in complex hyperbolic
space ${\Bbb C}H^n$ by Montiel and Romero \cite{MR}, and in compact complex two-plane Grassmannian $G_2({\Bbb C}^{m+2})$ by Berndt and Suh \cite{BS2} respectively. Moreover, for further investigating on commuting problems related to shape operator, Ricci tensor and the structure tensor are given. In complex projective space we want to mention some works due to Kimura \cite{K1}, \cite{K2}, in quaternionic projective space Martinez and P\'erez \cite{MP}, P\'erez and Suh \cite{PS1}, \cite{PS2}, and in complex two-plane Grassmannian $G_2({\Bbb C}^{m+2})$ P\'erez, Suh and Watanabe \cite{PSW} and  Suh \cite{S1}, \cite{S2} respectively.
\par
\vskip 6pt

The complex hyperbolic two-plane Grassmannian $\NBt$ has a remarkable
geometrical structure. It is the unique noncompact irreducible
Riemannian manifold being equipped with both a K\"{a}hler structure
$J$ and a quaternionic K\"{a}hler structure ${\frak J}=\
\text{Span}\ \{J_1, J_2, J_3\}$ not containing $J$. In other words,
$\NBt$ is the unique noncompact, irreducible, K\"{a}hler, quaternionic
K\"{a}hler manifold which is not a hyperk\"{a}hler manifold (See
Berndt and Suh \cite{BS3}).
\par
\vskip 6pt
\par
Now in this paper we want to give a characterization of type $(B)$, another one of type $(C)$, that is, a horosphere whose center at infinity is singular, or of type $(D)$ in noncompact complex two-plane Grassmannian $\NBt$ related to the Reeb vector field $\xi$. Then we can assert a complete classification of all Hopf real
hypersurfaces in  $\NBt$ in
terms of the {\it Reeb} vector field belonging to the maximal quaternionic
subbundle $\Cal Q$ as follows: \vskip 6pt

\par
\begin{Main Theorem}\label{Main Theorem}
Let $M$ be a Hopf real hypersurface in complex hyperbolic two-plane Grassmannian $\NBt$ with the
Reeb vector field belonging to the maximal quaternionic subbundle $\Cal Q$. Then one
of the following statements holds,
\par
$(B)$\quad $M$ is an open part of a tube around a totally geodesic
${\Bbb H}H^n$ in \hfill
\newline
$SU_{2,2n}/S(U_2U_{2n})$, $m = 2n$,
\par
$(C_2)$\quad $M$ is an open part of a horosphere in
$SU_{2,m}/S(U_2U_m)$ whose center at infinity is singular and of
type $JN \perp {\frak J}N$,
\par
or the following exceptional case holds:
\par
$(D)$\quad The normal bundle $\nu M$ of $M$ consists of singular
tangent vectors of type $JX \perp {\frak J}X$. Moreover, $M$ has at
least four distinct principal curvatures, three of which are given
by
\begin{equation*}
\alpha = \sqrt{2}\ ,\ \gamma = 0\ ,\ \lambda = \frac{1}{\sqrt{2}}
\end{equation*}
with corresponding principal curvature spaces
\begin{equation*}
T_\alpha = TM \ominus ({\Cal C} \cap {\Cal Q})\ ,\ T_\gamma = J(TM
\ominus {\Cal Q})\ ,\ T_\lambda \subset  {\Cal C} \cap {\Cal Q} \cap
J{\Cal Q}.
\end{equation*}
If $\mu$ is another (possibly nonconstant) principal curvature
function, then we have $T_\mu \subset {\Cal C} \cap {\Cal Q} \cap
J{\Cal Q}$, $JT_\mu \subset T_\lambda$ and ${\frak J}T_\mu \subset
T_\lambda $.
\end{Main Theorem}
\vskip 6pt
\par
\begin{rem}\quad Real hypersurfaces of type $(A)$ and $(C_1)$ with $JN{\in}{\frak J}N$ in
complex hyperbolic two-plane Grassmannians $\NBt$  mentioned in Theorem B
are characterized by the geometric property that the Reeb flow is isometric (See Suh \cite{S5}). Of course, in these type of hypersurfaces the Reeb vector field $\xi$ belongs
to the orthogonal complement ${\Cal Q}^{\bot}$ of the quaternionic maximal subbundle $\Cal Q$. But the other type of Hopf hypersurfaces in Theorem B are characterized in our main theorem
as Hopf hypersurfaces in $\NBt$ with the Reeb vector field ${\xi}{\in}{\Cal Q}$.
\end{rem}
\vskip 8pt
\par
\section{The complex hyperbolic two-plane
Grassmannian $\NBt$}\label{section 1}
\setcounter{equation}{0}
\renewcommand{\theequation}{1.\arabic{equation}}
\vspace{0.13in}

In this section we summarize basic material about the noncompact
complex two-plane Grassmann manifold $\NBt$, for details we refer to
\cite{BS3}, \cite{E}, \cite{H}, \cite{H2} and \cite{S5}.
\par
\vskip 6pt

The Riemannian symmetric space $SU_{2,m}/S(U_2U_m)$, which consists
of all positive definite complex two-dimensional linear subspaces in indefinite
complex Euclidean space $\CC_2^{m+2}$, becomes a connected, simply
connected, irreducible Riemannian symmetric space of noncompact type
with rank two. Let $G = SU_{2,m}$ and $K = S(U_2U_m)$, and
denote by ${\frak g}$ and ${\frak k}$ the corresponding Lie algebra.
Let $B$ be the Killing form of ${\frak g}$ and denote by ${\frak p}$
the orthogonal complement of ${\frak k}$ in ${\frak g}$ with respect
to $B$. The resulting decomposition ${\frak g} = {\frak k} \oplus
{\frak p}$ is a Cartan decomposition of ${\frak g}$. The Cartan
involution $\theta \in {\text Aut}({\frak g})$ on ${\frak s}{\frak
u}_{2,m}$ is given by $\theta(A) = I_{2,m} A I_{2,m}$, where
\begin{equation*}
I_{2,m} = \left(\begin{array}{cc}
 -I_2 & 0_{2,m} \\
 0_{m,2} & I_m
\end{array}\right)
\end{equation*}
and $I_2$ and $I_m$ is the identity $(2 \times 2)$-matrix and $(m \times m)$-matrix
respectively. Then $< X , Y > = -B(X,\theta Y)$ becomes a positive
definite ${\text Ad}(K)$-invariant inner product on ${\frak g}$. Its
restriction to ${\frak p}$ induces a Riemannian metric $g$ on
$SU_{2,m}/S(U_2U_m)$, which is also known as the Killing metric on
$SU_{2,m}/S(U_2U_m)$. Throughout this paper we consider
$SU_{2,m}/S(U_2U_m)$ together with this particular Riemannian metric
$g$.
\par
\vskip 6pt The Lie algebra ${\frak k}$ decomposes orthogonally into
${\frak k}  = {\frak s}{\frak u}_2 \oplus {\frak s}{\frak u}_m
\oplus {\frak u}_1$, where ${\frak u}_1$ is the one-dimensional
center of ${\frak k}$. The adjoint action of ${\frak s}{\frak u}_2$
on ${\frak p}$ induces the quaternionic K\"{a}hler structure ${\frak
J}$ on $SU_{2,m}/S(U_2U_m)$, and the adjoint action of
\begin{equation*}
Z = \left(\begin{array}{cc}
 \frac{mi}{m+2}I_2 & 0_{2,m} \\
 0_{m,2} & \frac{-2i}{m+2}I_m
 \end{array}\right)
 \in {\frak u}_1
\end{equation*}
induces the K\"{a}hler structure $J$ on $SU_{2,m}/S(U_2U_m)$. By
construction, $J$ commutes with each almost Hermitian structure
$J_1$ in ${\frak J}$. Recall that a canonical local basis
$J_1,J_2,J_3$ of a quaternionic K\"{a}hler structure ${\frak J}$
consists of three almost Hermitian structures $J_1,J_2,J_3$ in
${\frak J}$ such that $J_\nu J_{\nu+1} = J_{\nu + 2} = - J_{\nu+1}
J_\nu$, where the index $\nu$ is to be taken modulo $3$. The tensor
field $JJ_\nu$, which is locally defined on $SU_{2,m}/S(U_2U_m)$, is
selfadjoint and satisfies $(JJ_\nu)^2 = I$ and ${\text tr}(JJ_\nu) =
0$, where $I$ denotes the identity transformation. For a nonzero tangent
vector $X$ we define ${\Bbb R}X = \{\lambda X \vert \lambda \in
{\Bbb R}\}$, ${\Bbb C}X = {\Bbb R}X \oplus {\Bbb R}JX$, and ${\Bbb
H}X = {\Bbb R}X \oplus {\frak J}X$.
\par
\vskip 6pt Usually, the tangent space $T_oSU_{2,m}/S(U_2U_m)$ of
$SU_{2,m}/S(U_2U_m)$ at $o$ can be identified with ${\frak p}$. Let
${\frak a}$ be a maximal abelian subspace of ${\frak p}$. Since
$SU_{2,m}/S(U_2U_m)$ has rank two, the dimension of any such
subspace is two. Every nonzero tangent vector $X \in
T_oSU_{2,m}/S(U_2U_m) \cong {\frak p}$ is contained in some maximal
abelian subspace of ${\frak p}$. In general this subspace is
uniquely determined by $X$, in which case $X$ is called regular. If
there exists more than one maximal abelian subspaces of ${\frak p}$
containing $X$, then $X$ is called singular. There is a simple and
useful characterization of the singular tangent vectors: A nonzero
tangent vector $X \in {\frak p}$ is singular if and only if $JX \in
{\frak J}X$ or $JX \perp {\frak J}X$.
\par
\vskip 6pt Up to scaling there exists a unique $SU_{2,m}$-invariant Riemannian metric $g$ on $\NBt$. Equipped with this
metric $\NBt$ is a Riemannian symmetric space of rank two which is
both K\"ahler and quaternionic K\"ahler. For computational reasons
we normalize $g$ such that the minimal sectional curvature of
$({\NBt},g)$ is $-4$. The sectional curvature $K$ of the noncompact
symmetric space $SU_{2,m}/S(U_2U_m)$ equipped with the
Killing metric $g$ is bounded by $-4{\leq}K{\leq}0$. The sectional
curvature $-4$ is obtained for all $2$-planes ${\Bbb C}X$ when $X$
is a non-zero vector with $JX{\in}{\frak J}X$.
\par
\vskip 6pt

 When $m=1$,  $G_2^{*}(\CC^3)=SU_{1,2}/S(U_1U_2)$
is isometric to the two-dimensional complex hyperbolic space $\CC
H^2$ with constant holomorphic sectional curvature $-4$.
\par \vskip 6pt
When $m=2$, the isomorphism $SO(4,2)\simeq SU(2,2)$ yields an
isometry between $G_2^{*}(\CC ^4)=SU_{2,2}/S(U_2U_2)$ and the
indefinite real Grassmann manifold $G_2^{*}(\RR_2^6)$ of oriented
two-dimensional linear subspaces of an indefinite Euclidean space
$\RR_2^6$. For this reason we assume $m \geq 2$ from now on,
although many of the subsequent results also hold for $m = 1,2$.
\par
\vskip 6pt The Riemannian curvature tensor $\bar{R}$ of $\NBt$ is
locally given by
\begin{equation}\label{1.2}
\begin{split}
\bar{R}(X,Y)Z = &- \frac{1}{2}\Big[ g(Y,Z)X - g(X,Z)Y  +  g(JY,Z)JX\\
& - g(JX,Z)JY- 2g(JX,Y)JZ \\
& +  \sum_{\nu = 1}^3 \{g(J_\nu Y,Z)J_\nu X - g(J_\nu X,Z)J_{\nu}Y\\
& - 2g(J_\nu X,Y)J_\nu Z\} \\
& +  \sum_{\nu = 1}^3 \{g(J_\nu JY,Z)J_\nu JX - g(J_\nu JX,Z)J_{\nu}JY\}\Big ] ,
\end{split}
\end{equation}
where $J_1,J_2,J_3$ is any canonical local basis of ${\frak J}$.
\par
\vskip 6pt

Recall that a maximal flat in a Riemannian symmetric space $\bar{M}$
is a connected complete flat totally geodesic submanifold of maximal
dimension. A non-zero tangent vector $X$ of $\bar{M}$ is singular if
$X$ is tangent to more than one maximal flat in $\bar{M}$, otherwise
$X$ is regular. The singular tangent vectors of $\NBt$ are precisely
the eigenvectors and the asymptotic vectors of the self-adjoint
endomorphisms $JJ_1$, where $J_1$ is any almost Hermitian structure
in $\frak J$. In other words, a tangent vector $X$ to $\NBt$ is
singular if and only if $JX \in \frak J X$ or $JX \bot \frak J X$.
\par
\vskip 6pt In the previous paper of \cite{S5}, we considered a singular
vector of type $JX{\in}{\frak J}X$ and independently have given a
characterization of hypersurfaces of type $(A)$ and a horosphere of
type $(C_1)$. In this paper, we must compute explicitly Jacobi
vector fields along geodesics whose tangent vectors are all singular
of type $JX \perp \frak J X$. For this we need the eigenvalues and
eigenspaces of the Jacobi operator $\bar{R}_X := \bar{R}(.,X)X$. Let
$X$ be a singular unit vector tangent to $\NBt$ of type $JX \perp
\frak J X$.

 If $JX \perp {\frak J}X$ then the eigenvalues and eigenspaces of
$\bar{R}_X$ are given by(See Berndt and Suh \cite{BS3})
\begin{equation*}
\begin{array}{cll}
0 & \qquad \RR X \oplus {\frak J}JX & 4\\
-\frac{1}{2} & \qquad (\RR X {\oplus}\RR JX \oplus {\frak J}X \oplus {\frak J}JX)^{\bot} & 4m-8\\
-2 & \qquad \RR JX {\oplus} {\frak J}X & 4
\end{array}
\end{equation*}
where ${\RR}X$, ${\CC}X$ and ${\HH}X$ denote the real, complex and
quaternionic span of $X$, respectively, and $\CC ^\perp X$ the
orthogonal complement of $\CC X$ in $\HH X$. The maximal totally
geodesic submanifolds of $\NBt$ are $\NBo$, $\CC H^m$, $\CC H^k
\times \CC H^{m-k}$ ($1 \leq k \leq [m/2]$), $G_2^{*}(\RR^{m+2})$
and $\HH H^n$ (if $m = 2n$). The first three are complex
submanifolds and the other two are real submanifolds with respect to
the K\"ahler structure $J$. The tangent spaces of the totally
geodesic $\CC H^m$ are precisely the maximal linear subspaces of the
form $\{X \vert JX = J_1X\}$ with some fixed almost Hermitian
structure $J_1 \in \frak J$.

\vskip 8pt

\section{Real hypersurfaces in noncompact Grassmannian $\NBt$}\label{section 2}
\setcounter{equation}{0}
\renewcommand{\theequation}{2.\arabic{equation}}
\vspace{0.13in}

Let $M$ be a real hypersurface in $\NBt$, that is, a hypersurface in
$\NBt$ with real codimension one. The induced Riemannian metric on
$M$ will also be denoted by $g$, and $\nabla$ denotes the Levi
Civita covariant derivative of $(M,g)$. We denote by $\Cal C$ and
$\Cal Q$ the maximal complex and quaternionic subbundle of the
tangent bundle $TM$ of $M$, respectively. Now let us put
\begin{equation}\label{2.1}
JX={\phi}X+{\eta}(X)N,\quad J_{\nu}X={\phi}_{\nu}X+{\eta}_{\nu}(X)N
\end{equation}
for any tangent vector field $X$ of a real hypersurface $M$ in
$\NBt$, where ${\phi}X$ denotes the tangential component of $JX$ and
$N$ a unit normal vector field of $M$ in $\NBt$. From the K\"{a}hler
structure $J$ of $\NBt$ there exists an almost contact metric
structure $(\phi,\xi,\eta,g)$ induced on $M$ in such a way that
\begin{equation}\label{2.2}
{\phi}^2X=-X+{\eta}(X){\xi},\ {\eta}({\xi})=1,\ {\phi}{\xi}=0, \quad
\text{and}\quad {\eta}(X)=g(X,{\xi})
\end{equation}
for any vector field $X$ on $M$ and ${\xi}=-JN$.
\par
\vskip 6pt

If $M$ is orientable, then the vector field $\xi$ is globally
defined and said to be the induced {\it Reeb vector field} on $M$.
Furthermore, let $J_1,J_2,J_3$ be a canonical local basis of $\frak
J$. Then each $J_\nu$ induces a local almost contact metric
structure $(\phi_\nu,\xi_\nu,\eta_\nu,g)$, ${\nu}=1,2,3$, on $M$.
Locally, $\Cal C$ is the orthogonal complement in $TM$ of the real
span of $\xi$, and $\Cal Q$ the orthogonal complement in $TM$ of the
real span of $\{\xi_1,\xi_2,\xi_3\}$.
\par
\vskip 6pt

Furthermore, let $\{J_1,J_2,J_3\}$ be a canonical local basis of
${\frak J}$. Then the quaternionic K\"{a}hler structure $J_\nu$ of
$\NBt$, together with the condition
\begin{equation*}
J_{\nu}J_{\nu+1} = J_{\nu+2} = -J_{\nu+1}J_{\nu}
\end{equation*}
in section $1$, induced an almost contact metric 3-structure
$(\phi_{\nu}, \xi_{\nu}, \eta_{\nu}, g)$ on $M$ as follows:
\begin{equation}\label{2.3}
\begin{split}
&{\phi}_{\nu}^2X=-X+{\eta}_{\nu}({\xi}_{\nu}),\ {\phi}_{\nu}{\xi}_{\nu}=0,\ {\eta}_{\nu}({\xi}_{\nu})=1\\
&{\phi}_{\nu +1}{\xi}_{\nu}=-{\xi}_{{\nu}+2},\quad {\phi}_{\nu
}{\xi}_{{\nu}+1}={\xi}_{{\nu}+2},\\
&{\phi}_{\nu}{\phi}_{{\nu}+1}X
= {\phi}_{{\nu}+2}X+{\eta}_{{\nu}+1}(X){\xi}_{\nu},\\
&{\phi}_{{\nu}+1}{\phi}_{\nu}X=-{\phi}_{{\nu}+2}X+{\eta}_{\nu}(X)
{\xi}_{{\nu}+1}
\end{split}
\end{equation}
for any vector field $X$ tangent to $M$. The tangential and normal
component of the commuting identity $JJ_\nu X = J_\nu JX$ give
\begin{equation}\label{2.4}
\phi\phi_\nu X - \phi_\nu \phi X = \eta_\nu(X)\xi -
\eta(X)\xi_\nu \ \text{and}\quad \eta_\nu(\phi X) = \eta(\phi_\nu X).
\end{equation}

 The last equation implies
$\phi_\nu \xi = \phi \xi_\nu$. The tangential and normal component
of $J_\nu J_{\nu+1}X = J_{\nu+2}X = -J_{\nu+1}J_\nu X$ give
\begin{equation}\label{2.5}
\phi_\nu\phi_{\nu+1}X - \eta_{\nu+1}(X)\xi_\nu = \phi_{\nu+2}X
= - \phi_{\nu+1}\phi_\nu X + \eta_\nu(X) \xi_{\nu+1}
\end{equation}
and
\begin{equation}\label{2.6}
\eta_\nu(\phi_{\nu+1} X) = \eta_{\nu+2}(X) = - \eta_{\nu+1}(\phi_\nu X).
\end{equation}

Putting $X = \xi_\nu$ and $X = \xi_{\nu+1}$ into the first of these
two equations yields $\phi_{\nu+2}\xi_\nu = \xi_{\nu+1}$ and
$\phi_{\nu+2}\xi_{\nu+1} = - \xi_\nu$ respectively. Using the Gauss
and Weingarten formulas, the tangential and normal component of the
K\"ahler condition $(\bar{\nabla}_XJ)Y = 0$ give $(\nabla_X\phi)Y =
\eta(Y)AX - g(AX,Y)\xi$ and $(\nabla_X\eta)Y = g(\phi AX,Y)$. The
last equation implies $\nabla_X\xi = \phi AX$. Finally, using the
explicit expression for the Riemannian curvature tensor $\bar{R}$ of
$\NBt$ in (\ref{1.2}) the Codazzi equation takes the form
\begin{equation}\label{2.7}
\begin{split}
(\nabla_XA)Y &- (\nabla_YA)X
 = -\frac{1}{2}\Big[\eta(X)\phi Y - \eta(Y)\phi X - 2g(\phi X,Y)\xi \\
& \qquad +  \sum_{\nu = 1}^3 \big\{\eta_\nu(X)\phi_\nu Y - \eta_\nu(Y)\phi_\nu
 X - 2g(\phi_\nu X,Y)\xi_\nu\big\} \\
& \qquad +  \sum_{\nu = 1}^3 \big\{\eta_\nu(\phi X)\phi_\nu\phi Y
- \eta_\nu(\phi Y)\phi_\nu\phi X\big\} \\
& \qquad +  \sum_{\nu = 1}^3 \big\{\eta(X)\eta_\nu(\phi Y) - \eta(Y)\eta_\nu(\phi
X)\big\}\xi_\nu \Big] .
\end{split}
\end{equation}

\vskip 6pt We now assume that the Reeb flow on $M$ in $\NBt$ is
geodesic. Then, according to Proposition 3.1, there exists a smooth
function $\alpha$ on $M$ so that $A\xi = \alpha\xi$. Taking an inner
product of the Codazzi equation (\ref{2.7}) with $\xi$ we get
\begin{equation}\label{2.8}
\begin{split}
g(\phi X,Y)& - \sum_{\nu = 1}^3 \{\eta_\nu(X)\eta_\nu(\phi Y) - \eta_\nu(Y)\eta_\nu(\phi X)\\
& - g(\phi_\nu X,Y)\eta_\nu(\xi)\} \\
=& g((\nabla_XA)Y - (\nabla_YA)X, \xi)\\
=& g((\nabla_XA)\xi,Y) - g((\nabla_YA)\xi,X) \\
=&(X\alpha)\eta(Y) - (Y\alpha)\eta(X)\\
& + \alpha g((A\phi + \phi A)X,Y) - 2g(A\phi AX,Y)\ .
\end{split}
\end{equation}

Substituting $X = \xi$ yields $Y\alpha =
(\xi\alpha)\eta(Y) + 2\sum_{\nu = 1}^3 \eta_\nu(\xi)\eta_\nu(\phi
Y)$, and inserting this equation and the corresponding one for
$X\alpha$ into the previous equation implies
\vskip 6pt
\medskip
\begin{proposition}\label{Proposition 2.1}
If $M$ is a connected orientable
real hypersurface in \hfill
\newline
$\NBt$ with geodesic Reeb flow, then
\begin{equation*}
\begin{split}
 2&g(A\phi AX,Y) - \alpha g((A\phi + \phi A)X,Y) +  g(\phi X,Y)\\
 =&  \sum_{\nu=1}^3 \{\eta_\nu(X)\eta_\nu(\phi Y) - \eta_\nu(Y)\eta_\nu(\phi X) - g(\phi_\nu X,Y)\eta_\nu(\xi)\\
 & - 2 \eta(X)\eta_\nu(\phi Y) \eta_\nu(\xi) + 2 \eta(Y) \eta_\nu(\phi X) \eta_\nu(\xi)\}\ .
\end{split}
\end{equation*}
\end{proposition}
\vskip 6pt
\par

Here after, unless otherwise stated, we want to use these basic
equations mentioned above frequently without referring to them
explicitly.
\vskip 6pt
\par

\section{Proof of Main Theorem}\label{section 3}
\setcounter{equation}{0}
\renewcommand{\theequation}{3.\arabic{equation}}
\vspace{0.13in}

Let $M$ be a connected orientable Hopf real hypersurface in $\NBt$.
Now we denote by the distribution $\Cal Q$ the orthogonal
complement of the distribution
${\Cal Q}^{\bot}=\text{Span}\{\,\xi_{1}, \xi_{2}, \xi_{3}\,\}$ such
that $T_{x}M= \Cal Q \oplus {\Cal Q}^{\bot}$ for any point $x \in
M$.

In order to prove our Main Theorem in the introduction we give a key
proposition as follows:
\par
\vskip 6pt
\begin{proposition}\label{Proposition 3.1}
Let $M$ be a connected
orientable Hopf real hypersurface in \\
$\NBt$. If the Reeb vector field $\xi$ belongs to ${\Cal Q}$,
then $g(A {\Cal Q}, {\Cal Q}^{\bot})=0$.
\end{proposition}
\vskip 6pt

\begin{proof}\label{Proof}
To prove this it suffices to show that
$g(A\Cal{Q}, \xi_{\nu})=0,$ $\nu = 1, 2, 3$. In order to do this,
we put
\begin{equation*}
\Cal{Q}=
[\,\xi\,]\oplus[\,\phi_{1}\xi, \phi_{2}\xi, \phi_{3}\xi\,]\oplus
\Cal{Q}_{0},
\end{equation*}
where the distribution ${\Cal Q}_0$ is an
orthogonal complement of $[\,\xi\,]\oplus[\,\phi_{1}\xi,
\phi_{2}\xi, \phi_{3}\xi\,]$ in the distribution $\Cal Q$.
\par
First, from the assumption ${\xi}{\in}{\Cal Q}$ we know $g(A\xi,
\xi_{\nu})=0$, $\nu = 1,2,3,$ because we have assumed that $M$ is
Hopf.
\par
\vskip 6pt
Next we will show that $g(A\phi_{i} \xi, \xi_{\nu}) =0$,
for any indices $i$ and $\nu= 1,2,3$.
\par
In fact, by using (\ref{2.4}) and ${\xi}{\in}{\Cal Q}$ we have the
following
\begin{equation*}
\begin{split}
g(A \phi_{i} \xi , \xi_{\nu})&= -g(\phi A \xi_{\nu}, \xi_{i})\\
& = -g (\nabla_{\xi_{\nu}}\, \xi, \xi_{i}) \\
& = g( \xi, \nabla_{\xi_{\nu}}\,\xi_{i})\\
& = g( \xi, q_{i+2}(\xi_{\nu})\xi_{i+1}-q_{i+1}(\xi_{\nu})\xi_{i+2} + \phi_{i} A
\xi_{\nu})\\
& = g( \xi, \phi_{i}A \xi_{\nu}) \\
& =-g(A\phi_{i}\xi, \xi_{\nu}),
\end{split}
\end{equation*}
which gives that $g(A \phi_{i} \xi , \xi_{\nu}) = 0$, $\nu=1, 2, 3$.
\par
\vskip 6pt
Finally, we consider the case $X \in \Cal{Q}_{0}$, where the
distribution ${\Cal Q}_0$ is denoted by
\begin{equation*}
\Cal{Q}_{0} =\{X \in \Cal{Q}\,|\ X \bot\ \xi\ \text{and} \
\phi_{i}\xi\,, \ i=1,2,3 \}.
\end{equation*}
By Proposition 2.1 and the assumption of ${\xi}{\in}{\Cal Q}$,
we have
\begin{equation*}
\begin{split}
& \alpha A \phi X  +  \alpha \phi AX -2 A \phi AX  - \phi X \\
 =& -2 \sum_{\nu=1}^{3}\Big \{{\eta}(X){\ENK}{\PKN} + \eta_{\nu}(\phi X){\ENK}{\xi}\Big \}\\
&+ \sum_{\nu=1}^3\Big\{\eta_{\nu}(X)\phi_{\nu}{\xi}+\eta_{\nu}({\phi}X)\xi_{\nu}+\eta_{\nu}(\xi){\phi}_{\nu}X\Big\}\\
=&0
\end{split}
\end{equation*}
for any tangent vector field $X \in {\Cal Q}_0$.
\par
\vskip 6pt
From now on, in order to show $g(AX, \xi_{\nu})=0$ for any $X \in
\Cal{Q}_{0}$,  we restrict $X \in T_{p}M$, $p{\in}M$ to $ X \in \Cal{Q}_{0}$
unless otherwise stated. Now by taking the structure tensor $\phi$ into above equation
and using the fact that $\xi \in \Cal{Q}$ we get
\begin{equation}\label{3.1}
\begin{split}
\alpha \phi A \phi X & - \alpha AX -2 \phi A \phi AX  +  X =0 \,,
\end{split}
\end{equation}
for any $X \in \Cal{Q}_{0}$.

Taking an inner product into (\ref{3.1}) with $\xi_{\mu}$ we have
\begin{equation*}
\begin{split}
\alpha g(\phi A \phi X, \xi_{\mu})   - \alpha g(AX, \xi_{\mu}) -2 g(\phi A \phi AX, \xi_{\mu})=0 \,,
\end{split}
\end{equation*}
that is,
\begin{equation}\label{3.2}
\begin{split}
\alpha g(AX, \xi_{\mu})= \alpha g(\phi A \phi X, \xi_{\mu})-2 g(\phi A \phi AX, \xi_{\mu}) \quad \quad \text {for}\ X \in \frak{Q}_{0}.
\end{split}
\end{equation}

On the other hand, since $g(\phi A \phi X, \xi_{\mu})= g(\nabla_{\phi X}\xi , \xi_{\mu})= - g( \xi, \nabla_{\phi X} \xi_{\mu})$, we have
\begin{equation*}
g(\phi A \phi X, \xi_{\mu})= -g(\xi, \phi_{\mu} A \phi X)= -g(\xi_{\mu}, \phi A \phi X)
\end{equation*}
by virtue of (\ref{2.1}) and (\ref{2.4}). Accordingly, we get $g(\phi A \phi X,
\xi_{\mu})= 0$ for any $X \in \Cal{Q}_{0}$.
\par
\vskip 6pt
Next let us show that $g(\phi A \phi AX, \xi_{\mu})=0$.
\vskip 6pt
\par
In fact, (\ref{2.3}) and (\ref{2.4}) give
\begin{equation*}
\begin{split}
g(\phi A \phi AX, \xi_{\mu})=& g(\nabla_{\phi AX} \xi,
\xi_{\mu})= -g( \xi, \nabla_{\phi AX}\xi_{\mu})\\
=& -g( \xi, \phi_{\mu} A \phi AX )=-g(\xi_{\mu}, \phi A \phi AX).
\end{split}
\end{equation*}
It implies that $g(\phi A \phi AX, \xi_{\mu})=0$ for any $X \in
\Cal{Q}_{0}$. Thus, from (\ref{3.2}) we know that
\begin{equation}\label{3.3}
\alpha g(AX, \xi_{\mu})=0 \quad \text{ for \ any \ } X \in \Cal{Q}_{0}\,.
\end{equation}
\par
From this we can divide two cases as follow:
\par
\vskip 6pt
{\bf Case I.} Let ${\frak U}=\{x{\in}M{\vert}{\alpha}(x){\not =}0\}$.
\par
\vskip 6pt
On such an open neighborhood $\frak U$ we know that (\ref{3.3}) gives $g(AX, \xi_{\mu})=0$ for any $X
\in \Cal{Q}_{0}$.
\par
\vskip 6pt
{\bf Case II.} Let ${\frak W}=\text{Int}\ (M - {\frak U})$, where {\it Int} denotes the interior set of the orthogonal
complement of the open subset $\frak U$ in $M$.
\par
\vskip 6pt
In this case we consider two subcases. One is to consider that the fucntion $\alpha$ vanishes on a non-empty neighborhood $\text{Int}\ (M-{\frak U})$.
The other subcase is to consider a point $x$ such that ${\alpha}(x)=0$ but the point $x$ is the limit of a sequence of points where
${\alpha}{\not =}0$. This subcase could be possible when the open subset $\text{Int}\ (M-{\frak U})$ is empty. Such a sequence necessary have an infinite subsequence. Then by the continuity we have $g(AX, \xi_{\mu})=0$ for any $X
\in \Cal{Q}_{0}$ as in Case I.
\par
\vskip 6pt
Then we only focus on the first subcase that the the function ${\alpha}$ identically vanishes on some neighborhood of the point $x{\in}M$.
From this situation, the equation (\ref{3.1}) can be given by
\begin{equation}\label{3.4}
X= 2\phi A \phi AX  \quad \text{ for \ any \ } X \in \Cal{Q}_{0}\,.
\end{equation}

Taking the shape operator $A$ into (\ref{3.4}) we have
\begin{equation}\label{3.5}
AX= 2A\phi A \phi AX \quad \text{ for \ any \ } X \in \Cal{Q}_{0}\,.
\end{equation}

From this, let us take an inner product into (\ref{3.5}) with $\xi_{\mu}$,
we have
\begin{equation}\label{3.6}
g(AX, \xi_{\mu})= 2g( A\phi A \phi AX, \xi_{\mu}) \quad \text{ for \
any \ } X \in \Cal{Q}_{0}.
\end{equation}

On the other hand, we know the following
\begin{equation*}
g(A \phi A \phi AX,
\xi_{\mu}) = -g(A \phi A X, \phi A \xi_{\mu})= -g(A \phi AX,
\nabla_{\xi_{\mu}} \xi)
\end{equation*}
Then it follows that
\begin{equation*}
\begin{split}
g(A \phi A \phi AX, \xi_{\mu}) & = -g(A \phi A X, \nabla_{\xi_{\mu}} \xi)\\
=& g((\nabla_{\xi_{\mu}}A) \phi AX, \xi) + g(A(\nabla_{\xi_{\mu}}\phi)AX, \xi)\\
& + g(A \phi (\nabla_{\xi_{\mu}}A)X, \xi) + g(A \phi A (\nabla_{\xi_{\mu}}X), \xi)\\
=&g(({\nabla}_{\xi_{\mu}}A){\phi}AX,{\xi})
\end{split}
\end{equation*}
where we have used $g(A \phi AX, \xi)=0$ and $A{\xi}=0$. From this, together with
$A\xi = 0$, it follows that
\begin{equation}\label{3.7}
g(A \phi A \phi AX, \xi_{\mu}) = g((\nabla_{\xi_{\mu}}A) \phi AX, \xi)\,.
\end{equation}

\par
On the other hand, by using the equation of Codazzi in section 2, we
have the following
\par
\vskip 6pt
\begin{lemma}\label{Lemma 3.2}
\begin{equation*}
g((\nabla_{\xi_{\mu}}A)
\phi AX, \xi) =  -g(A\xi_{\mu}, \phi A \phi A X)+2g(AX,{\xi}_{\mu}), \quad {\mu}=1,2,3 .
\end{equation*}
\end{lemma}
\vskip 6pt

\begin{proof}\label{Proof}
 By using the equation of Codazzi, it follows that for $\xi \in \Cal{Q}$
\begin{equation*}
\begin{split}
(\nabla_{\xi_{\mu}}A) \phi AX &=  (\nabla_{\phi AX}A) \xi_{\mu} -\frac{1}{2}\Bigg[ \eta(\xi_{\mu})\phi^{2}AX - \eta(\phi AX) \phi \xi_{\mu} -2g( \phi \xi_{\mu}, \phi AX) \xi \\
&\quad + \sum_{\nu=1}^{3}\Big\{ \eta_{\nu} (\xi_{\mu}) \phi_{\nu} \phi AX - \eta_{\nu} (\phi AX) \phi_{\nu} \xi_{\mu} -2g(\phi_{\nu} \xi_{\mu}, \phi AX) \xi_{\nu}\Big \}\\
&\quad  + \sum_{\nu=1}^{3}\Big \{ \eta_{\nu}(\phi \xi_{\mu}) \phi_{\nu} \phi^{2}AX - \eta_{\nu}(\phi^{2}AX) \phi_{\nu}\phi \xi_{\mu}\Big \}\\
&\quad  +\sum_{\nu=1}^{3} \Big\{ \eta(\xi_{\mu}) \eta_{\nu}(\phi^{2}AX) - \eta(\phi AX) \eta_{\nu} (\phi \xi_{\mu})\Big \} \xi_{\nu}\Bigg]\\
& = (\nabla_{\phi AX}A) \xi_{\mu} + g(\xi_{\mu}, AX) \xi -\frac{1}{2}\phi_{\mu} \phi AX \\
& + \frac{1}{2}\sum_{\nu=1}^{3}{\eta}_{\nu}({\phi}AX)\phi_{\nu}\xi_{\mu}\\
& + \sum_{\nu=1}^{3} g(\phi_{\nu}\xi_{\mu}, {\phi}AX) \xi_{\nu}
-\frac{1}{2}\sum_{\nu=1}^{3} \eta_{\nu}(AX) \phi_{\nu} \phi \xi_{\mu}.
\end{split}
\end{equation*}
Taking an inner product above equation with $\xi$ and using the fact
that $\phi \phi_{\mu} \xi = -\xi_{\mu}$, we have the following for any $X{\in}{\Cal Q}_0$
\begin{equation*}
\begin{split}
g((\nabla_{\xi_{\mu}}A) \phi AX, \xi) =&  g((\nabla_{\phi AX}A) \xi_{\mu}, \xi) +g({\xi}_{\mu},AX) - \frac{1}{2}g({\phi}_{\mu}{\phi}AX,{\xi})   \\
& \quad   + \frac{1}{2}\sum_{\nu=1}^{3}{\eta}_{\nu}({\phi}AX)g(\phi_{\nu}\xi_{\mu}, \xi) + \sum_{\nu=1}^{3} g(\phi_{\nu}\xi_{\mu}, {\phi}AX) g(\xi_{\nu},\xi)\\
& \quad - \frac{1}{2}\sum_{\nu=1}^{3} \eta_{\nu}(AX)g(\phi_{\nu} \phi \xi_{\mu},\xi)\\
=& g((\nabla_{\phi A X} A) \xi_{\mu}, \xi) + \frac{3}{2}g(AX, \xi_{\mu}) + \frac{1}{2}{\eta}_{\mu}(AX)\\
=& g((\nabla_{\phi A X} A) \xi_{\mu}, \xi) + 2g(AX,{\xi}_{\mu}),
\end{split}
\end{equation*}
where we have used $g(\phi_{\nu}\xi_{\mu}, \xi)=0$ and $g(\xi_{\nu},\xi)=0$ in the second equality.
\par
\vskip 6pt
On the other hand, since $g(A \xi_{\mu}, \xi) = g(\xi_{\mu}, A\xi) =
\alpha g (\xi_{\mu}, \xi)$ and $\alpha =0$, we have
\begin{equation*}
\begin{split}
g((\nabla_{\phi A X} A) \xi_{\mu}, \xi)&= -g(A (\nabla_{\phi AX}\xi_{\mu}),\xi) - g(A \xi_{\mu}, \phi A \phi AX) \\
& = - \alpha g( \nabla_{\phi A X }\xi_{\mu}, \xi) - g(A\xi_{\mu}, \phi A \phi AX)\\
& = - g(A\xi_{\mu}, \phi A \phi AX) \,.
\end{split}
\end{equation*}
Therefore we have
\begin{equation*}
g((\nabla_{\xi_{\mu}}A) \phi AX, \xi)= - g(A\xi_{\mu}, \phi A \phi AX) + 2g(AX,{\xi}_{\mu})
\end{equation*}
for any $X \in \Cal{Q}_{0}$. This completes the proof of our Lemma 3.2.
\end{proof}
\par

\par
Consequently, from (\ref{3.7}), Lemma 3.2 and the formula $g(A\xi_{\mu},
\xi)=0$ we have
\begin{equation*}
\begin{split}
g(A \phi A \phi AX, \xi_{\mu})& = g((\nabla_{\xi_{\mu}}A) \phi AX, \xi)\\
&= -g(A\xi_{\mu} , {\phi}A{\phi}AX)) + 2g(AX,{\xi}_{\mu}) \\
& = -g(A \phi A \phi AX, \xi_{\mu}) + 2g(AX,{\xi}_{\mu})
\end{split}
\end{equation*}
that is,
\begin{equation}\label{3.8}
g(A \phi A \phi AX, \xi_{\mu}) = g(AX,{\xi}_{\mu}).
\end{equation}

Summing up (\ref{3.6}) and (\ref{3.8}) for ${\alpha}=0$, we have $g(AX,
\xi_{\mu})=0$. Then for any $X \in \Cal{Q}_{0}$ we have $g(AX,
\xi_{\mu}) = 0$, $\mu= 1, 2, 3\,$. This completes the proof of our
Proposition 3.1.
\end{proof}
\vskip 6pt
\par
 By virtue of Proposition 3.1, the maximal complex subbundle $\Cal C$ and the maximal quaternionic subbundle $\Cal Q$ of Hopf real hypersurfaces $M$ in $\NBt$ are invarinat under the shape opeartor if the Reeb vector field ${\xi}$ of $M$ belongs to the subbundle ${\Cal Q}$. Then naturally we get the result of Theorem B. But among the classifications
given in Theorem B, the tube over a totally geodesic
$SU_{2,m-1}/S(U_2U_{m-1})$ in $SU_{2,m}/S(U_2U_m)$ and  a horosphere in $SU_{2,m}/S(U_2U_m)$
 whose center at infinity with singular vector field of type $JX{\in}{\Cal J}X$ have the property that their Reeb vector field ${\xi}$ belong to the subbundle ${\Cal Q}^{\bot}$ which is orthogonal to the maximal quaternionic subbundle $\Cal Q$. From such a point of view we complete
 the proof of our Main Theorem in the introduction.
\par
\vskip 8pt
{\bf Acknowledgements}\quad The present author would like to express his deep gratitude to Professor J\"urgen Berndt for his valuable comments and suggestions to develop the first version of this manuscript.



\begin{thebibliography}{99}

\bibitem{A} D.V. Alekseevskii, {\it Compact quaternion spaces}, Func. Anal. Appl. ~{\bf 2} (1968), 106--114.

\bibitem{Be} A. L. Besse, {\it Einstein manifolds}, Classics in Mathematics, Reprint of the 1987 Edition with 22 Fgures, Springer-Verlag, \textbf{2008}.


\bibitem{BS2} J. Berndt and Y.J. Suh, {\it Real hypersurfaces with isometric Reeb flow in complex two-plane Grassmannians}, Monatshefte f\"ur Math. ~{\bf 137} (2002), 87--98.

\bibitem{BS3} J. Berndt and Y.J. Suh, {\it Hypersurfaces in noncompact complex Grassmannians of rank two}, International J. Math., World Scientific Publishing, {\bf 23}(2012), 1250103(35 pages).

\bibitem{E} P.B. Eberlein, {\it Geometry of nonpositively curved manifolds}, University of Chicago Press, Chicago, London. \textbf{1996}.

\bibitem{H} S. Helgason, {\it Groups and Geometric Analysis}, Survey and Monographs Amer. Math. Soc. ~{\bf 83} 2002.

\bibitem{H2} S. Helgason, {\it Geometric Analysis on Symmetric Spaces}, The 2nd Edition, Math. Survey and Monographs, Amer. Math. Soc. ~{\bf 39} (2008).

\bibitem{K1} M.\ Kimura,  {\it Real hypersurfaces and complex submanifolds in complex projective space}, Trans. Amer. Math. Soc. ~{\bf 296} (1986), 137--149.

\bibitem{K2} M. Kimura, {\it Some real hypersurfaces of a complex projective space}, Saitama Math. J. {\bf 5} (1987), 1-5.


\bibitem{MP}  A. Martinez and J.D. P{\'e}rez, {\it Real hypersurfaces in quaternionic projective space}, Ann. di Mate. Pura Appl. ~{\bf 145} (1986), 355--384.

\bibitem{MR} S. Montiel and A. Romero, {\it On some real hypersurfaces of a complex hyperbolic space}, Geometriae Dedicata ~{\bf 20} (1986), 245--261.

\bibitem{O} M. Okumura, {\it On some real hypersurfaces of a complex projective space}, Trans. Amer. Math. Soc. ~{\bf 212} (1975), 355--364.



\bibitem{PS1}  J.D. P\'erez and Y.J. Suh, {\it Real hypersurfaces of
quaternionic projective space satisfying ${\nabla}_{U_i}R=0$}, Diff. Geom. and Its Appl. {\bf 7} (1997), 211--217.
\noindent

\bibitem{PS2}  J.D. P\'erez and Y.J. Suh, {\it Certain conditions on the Ricci tensor of real hypersurfaces in quaternionic projective space}, Acta Math. Hungarica {\bf 91} (2001), 343-356.
\noindent

\bibitem{PSW} J.D. P\'erez, Y.J. Suh and Y. Watanabe, {\it Generalized Einstein real hypersurfaces in
complex two-plane Grassmannians}, J. of Geom. and Physics
\textbf{60} (2010), 1806-1818.
\noindent%

\bibitem{S1}  Y.J. Suh, {\it Real hypersurfaces of type $B$ in
complex two-plane Grassmannians}, Monatshefte f\"ur Math. {\bf 147} (2006), 337--355.
\noindent

\bibitem{S2}  Y.J. Suh, {\it Real hypersurfaces in complex two-plane Grassmannians with commuting Ricci tensor}, J. Geometry and Physics {\bf 60} (2010), 1792-1805.
\noindent




\bibitem{S5} Y. J. Suh, {\it Real hypersurfaces with isometric Reeb flow in complex hyperbolic two-plane Grassmannians},  Adv. in Applied Math. {\bf 50} (2013), 645-659.

\end{thebibliography}
\end{document}